\documentclass[a4paper]{article}
\usepackage{amsmath,amsfonts,amssymb,amsthm,amstext}
\usepackage{mathtools}
\usepackage{latexsym}
\usepackage{microtype}
\usepackage{lmodern}
\usepackage[T1]{fontenc}
\usepackage[utf8]{inputenc}
\usepackage[english]{babel}
\xdef\articletitle{Parallel Type Decomposition Scheme for Quasi-Linear Abstract Hyperbolic Equation}
\xdef\Dikhaminjia{Nana Dikhaminjia}
\xdef\Rogava{Jemal Rogava}
\xdef\Tsiklauri{Mikheil Tsiklauri}
\usepackage[pdftex,
pdfauthor={ {\Dikhaminjia}, {\Rogava} and {\Tsiklauri} },
pdftitle={\articletitle},
pdfsubject={Mathematics, Applied Mathematics},
pdfstartview={FitH}]{hyperref}
\usepackage{mathscinet}
\usepackage{enumitem}
\usepackage{graphicx}
\usepackage{float}
\usepackage{color}
\usepackage[nottoc]{tocbibind}
\usepackage[square,sort,comma,numbers,sectionbib]{natbib}
\graphicspath{ {Figures/} }
\usepackage{geometry}
\numberwithin{equation}{section}

\theoremstyle{plain}

\newtheorem*{definition*}{Definition}



\providecommand{\keywords}[1]
{
	\small	
	\noindent\textbf{\textit{Keywords and phrases:}} #1
}

\providecommand{\msc}[1]
{
	\small	
	\noindent\textbf{\textit{MSC 2010:}} #1
}

\title{\articletitle}
\author{ {\Dikhaminjia}, {\Rogava} and {\Tsiklauri} }

\date{}

\begin{document}

\maketitle

\begin{abstract}
Cauchy problem for an abstract hyperbolic equation with the Lipschitz continuous operator is considered in the Hilbert space. The operator corresponding to the elliptic part of the equation is a sum of operators $A_{1},\,A_{2},\,\ldots,\,A_{m}$. Each addend is a self-adjoint and positive definite operator. A parallel type decomposition scheme for an approximate solution of the stated problem is constructed. The main idea of the scheme is that on each local interval classic difference problems are solved in parallel (independently from each other) respectively with the operators $A_{1},\,A_{2},\,\ldots,\,A_{m}$. The weighted average of the received solutions is announced as an approximate
solution at the right end of the local interval. Convergence of the proposed scheme is proved and the approximate solution error is estimated, as well as the error of the difference analogue for the first-order derivative for the case when the initial problem data satisfy the natural sufficient conditions for solution existence.
\end{abstract}
\keywords{Decomposition scheme, Abstract hyperbolic equation, Operator splitting, Cauchy problem, Parallel algorithm.} \\
\msc{65M12, 65M15, 65M55, 49M27.}

\section{Introduction}

\label{S1} \vspace{-4pt}

First works dedicated to construction and investigation of decomposition
schemes were published in the fifties and sixties of the XX century (see G. A. Baker, T. A. Oliphant \cite{BO}, G. Birkhoff, R. S. Varga \cite{BirV}, G. Birkhoff, R. S. Varga, D. Young \cite{BirVY}, J. Douglas \cite{Do}, J.
Douglas, H. Rachford \cite{DoRa}, E. G. Diakonov \cite{Dia}, D. G. Gordeziani \cite{DG}, N. N. Ianenko \cite{Ia}, \cite{Ia2}, V. P. Ilin \cite{Ilin}, A. N. Konovalov \cite{Ko}, G. I. Marchuk, N. N. Ianenko \cite{MaIa},
G. I. Marchuk, U. M. Sultangazin \cite{MaSu}, D. Peaceman, H. Rachford \cite{PeRa}, A. A. Samarskii \cite{Sa}, \cite{Sa2}). It may be said that the works of these authors become a basis for further research on decomposition schemes.

Decomposition schemes in view of numerical calculation can be divided into two groups: schemes of sequential account (see for example G. I. Marchuk \cite{Mar}, A. A. Samarskii, P. N. Vabishchevich \cite{SV} ) and schemes of parallel account (D. G. Gordeziani \cite{DG}, \cite{DG2}, D. G. Gordeziani, H. V. Meladze \cite{GM}, D. G. Gordeziani, A. A. Samarskii \cite{DGS}, A. M. Kuzyk, V. L. Makarov \cite{KMa}).

In view of numerical calculations, with the development of parallel processing, obviously the parallel type decomposition schemes have clear advantage. The present work discusses the construction and investigation of parallel type decomposition scheme.

It is well-known that for error estimation of approximate solution of evolution problem, usually the solution is required to be of the higher order smooth than it is necessary following from natural conditions. In case of decomposition schemes this detail gains more importance. Demand on increasing the smoothness can be conditioned by the fact that the operators corresponding to the problems obtained by splitting are non-commutative. Therefore, it is important to build the decomposition schemes, whose numerical calculation and methodology of error estimate of approximate solution does not require sharp increase of the solution smoothness. These details are much more complicated for the second order evolution equation compared to the first
order one. One of the reasons for this is that natural scheme for the first order evolution equation is two-layer, and for the second order - three-layer. For the most cases, in comparison with the two-layer schemes, investigation of the three-layer schemes are related with certain difficulties. At the first sight, this issue might be overcome: the second order evolution equation by introducing the additional unknown can be deduced to the first order system. However, in this case, if the operator in the initially given equation is self-adjoint, in the obtained system there will be non self-adjoint matrix operator, that significantly complicates solving of the corresponding discrete problem.

In the present work the decomposition scheme for the second order evolution equation is proposed that does not require to increase smoothness of the solution in view of numerical calculation. In addition, the methodology that we use for error estimate of the approximate solution makes it possible to find convergence order in the conditions with almost natural limitations.

For the investigation of the decomposition scheme we
use polynomials of certain class, which we call two-variable polynomial. These polynomials are represented by means of second order classical Chebyshev polynomials.

We should note that several works are devoted to use of orthogonal polynomials in approximate solution schemes for differential equations: V. L. Makarov \cite{Mak}, A. G. Morris, T. S. Horner \cite{MoHo}, V. A. Novikov, G. V. Demidov \cite{ND}, V. A. Rastrenin \cite{Ras}. In the work \cite{Mak} many aspects of using orthogonal polynomials
in the difference problems are presented quite widely.

Recent results related to the construction and investigation of decomposition schemes for evolution equations are obtained by the following authors: S. Blanes, F. Casas and M. Thalhammer \cite{Sergio}, D. He, K. Pan and H. Hu \cite{Dongdong}, J. L. Padgett and Q. Sheng \cite{Joshua}, J. Zhao, R. Zhan and Y. Xu \cite{Jingjun}. We also note the work \cite{DikhRogaTsik2014} in which a high-order accuracy decomposition scheme is considered for an abstract hyperbolic equation.

\section{Statement of the problem and decomposition scheme}

Let us consider the Cauchy problem for abstract hyperbolic equation in the
Hilbert space $H$:

\begin{equation}
\frac{d^{2}u(t)}{dt^{2}}+Au\left( t\right) +M\left( u\left( t\right) \right)
=f\left( t\right) ,\quad t\in \left[ 0,T\right] ,  \label{2.1}
\end{equation}%
\begin{equation}
u\left( 0\right) =\varphi _{0},\quad \frac{du\left( 0\right) }{dt}=\varphi
_{1}.  \label{2.1_2}
\end{equation}%
where $A$ is a self-adjoint ($A$ does not depend on $t$), positive definite
(generally unbounded) operator with the definition domain $D\left( A\right) $%
, which is everywhere dense in $H$, i.e. $\overline{D\left( A\right) }%
=H,\quad A=A^{\ast }$ and
\[
\left( Au,u\right) \geq \alpha \left\Vert u\right\Vert ^{2},\quad \forall
u\in D\left( A\right) ,\quad \alpha =const>0,
\]%
where by $\left\Vert \cdot \right\Vert $ and $\left( \cdot ,\cdot \right) $
are defined correspondingly the norm and scalar product in $H$; nonlinear
operator $M\left( \cdot \right) $ satisfies Liptschitz condition,

\[
\left\Vert M(u)-M(v)\right\Vert \leq a\left\Vert u-v\right\Vert ,\quad
\forall u,v\in H,
\]%
$a=const>0$ ; $\varphi _{0}$ and $\varphi _{1}$ are given vectors from $H$; $%
u\left( t\right) $ is a continuous, twice continuously differentiable,
searched function with values in $H$, and $f\left( t\right) $ is given
continuous function with values in $H.$

Similar to the linear case, $u\left( t\right) $ vector function with values
in $H$, defined on the interval $\left[ 0,T\right] $, is called a solution
of the problem (\ref{2.1}), (\ref{2.1_2}) if it satisfies the following
conditions: (a) $u\left( t\right) $ is twice continuously differentiable in
the interval $\left[ 0,T\right] $; (b) $u\left( t\right) \in D\left(
A\right) $ for any $t$ from $\left[ 0,T\right] $ , the function $\ Au\left(
t\right) $ is continuous and $M(u\left( t\right) )$ is continuous; (c) $%
u\left( t\right) $ satisfies equation (\ref{2.1}) on the $\left[ 0,T\right] $
interval and the initial condition (\ref{2.1_2}).

\textbf{Remark 2.1. }If $f\left( t\right) $ is continuously differentiable
on $\left[ 0,T\right] $ (or $f\left( t\right) \in D(A^{1/2})$ for any $t$
from $\left[ 0,T\right] $ and the function $A^{1/2}f\left( t\right) $ is
continuous), $\varphi _{0}\in D(A)$ and $\varphi _{1}\in D(A^{1/2})$, \ then
there exists only solution $u\left( t\right) $ of the problem (\ref{2.1}), (%
\ref{2.1_2}) (without Lipschitz continuous operator) that satisfies the
condition: the function $u^{\prime }\left( t\right) $ gets the values from $%
D(A^{1/2})$ and $A^{1/2}u^{\prime }\left( t\right) $ is continuous on $\left[
0,T\right] $ (see \cite{Kr}, Theorem 1.5, p. 301).

Let
\begin{equation}
A=\sum_{j=1}^{m}A_{j}\ ,\quad A_{j}=A_{j}^{\ast }\geq \alpha _{j}I\ ,\quad
\alpha _{j}=const>0\ .  \label{A.1}
\end{equation}

\qquad Then approximate solution of problem (\ref{2.1}), (\ref{2.1_2}) at
the points $t=t_{k+1}=(k+1)\tau $\ , $k=1,\ldots ,n-1$\ , \ $\tau =T/n$ ($%
n>1 $) is defined by the following formula:

\[
v_{k+1}=\sum_{j=1}^{m}\eta _{j}y_{j,k+1}\ ,\quad \sum_{j=1}^{m}\eta
_{j}=1,\quad 0<\eta _{j}<1,\quad
\]

where $y_{j,k+1}$ is a solution of the following difference problem:

\begin{equation}
\eta _{j}\frac{y_{j,k+1}-2v_{k}+v_{k-1}}{\tau ^{2}}+A_{j}y_{j,k+1}=\delta
_{1,j}\left[ f(t_{k})-M\left( v_{k}\right) \right] \ ,  \label{E1.1}
\end{equation}

\begin{equation}
v_{0}=\varphi _{0}\ ,\quad v_{1}=\varphi _{0}+\tau \varphi _{1}\ ,
\label{E1.2}
\end{equation}

where $j=1,\ldots ,m$, $\delta _{1,j}$ is a Kronecker symbol.

\qquad Thus, to construct approximate solution \ $v_{k+1}$ for problem (\ref%
{2.1}), (\ref{2.1_2}) at the point $t_{k+1}$ , it is necessary to solve $m$
\ problems independent from each other. Therefore, scheme (\ref{E1.1}) can
be called parallel type decomposition scheme. These kind schemes for the
first time were discussed in the work by D. Gordeziani (see \cite{DG}, \cite%
{DG2}). Parallel type decomposition schemes also are considered in the
works: D. Gordezianis, A. Samarski \cite{DGS}, D. G. Gordeziani, H. V.
Meladze \cite{GM}, A. M. Kuzyk, V. L. Makarov \cite{KMa}. Specifically the
scheme (\ref{E1.1}) (without Lipschitz-continuous operator) is given in \cite%
{RJ}.

\section{Representation of the approximate solution error by means of
Chebyshev polynomial and the main theorem}

Let the problem (\ref{2.1}), (\ref{2.1_2}) has a solution.

Then the equation (\ref{2.1}) at the point $t=t_{k+1}$ can be written as

\begin{equation}
\frac{u(t_{k+1})-2u(t_{k})+u(t_{k-1})}{\tau ^{2}}+Au(t_{k+1})=g_{k}\ ,
\label{E1.5}
\end{equation}
where
\begin{eqnarray*}
g_{k} &=&\widetilde{f}(t_{k})+A\left[ u(t_{k+1})-u(t_{k})\right] \\
&&+\tau ^{-2}\int_{t_{k}}^{t_{k+1}}(t_{k+1}-t)\left[ u^{\prime \prime
}\left( t\right) -u^{\prime \prime }\left( t_{k}\right) \right] dt
\end{eqnarray*}

\[
+\tau ^{-2}\int_{t_{k-1}}^{t_{k}}(t-t_{k-1})\left[ u^{\prime \prime }\left(
t\right) -u^{\prime \prime }\left( t_{k}\right) \right] dt\ .
\]
and where $\widetilde{f}(t)=f(t)-M(u(t))$.

From (\ref{E1.5}) it follows that
\begin{equation}
u(t_{k+1})-2Lu(t_{k})+Lu(t_{k-1})=\tau ^{2}Lg_{k}\ ,  \label{E1.7}
\end{equation}
where $k=1,\ldots ,n-1$,
\[
L=\left( I+\tau ^{2}A\right) ^{-1}.
\]

From (\ref{E1.1}) we have
\begin{equation}
y_{j,k+1}-2S_{j}v_{k}+S_{j}v_{k-1}=\tau ^{2}\eta _{j}^{-1}\delta _{1,j}S_{j}
\left[ f(t_{k})-M\left( v_{k}\right) \right] \ ,  \label{E1.8}
\end{equation}
where $j=1,\ldots ,m$ ,
\[
S_{j}=\left( I+\tau ^{2}\eta _{j}^{-1}A_{j}\right) ^{-1}.
\]

If we multiply both sides of equality (\ref{E1.8}) on $\eta _{j}$ and
summarize, we get
\begin{equation}
v_{k+1}-2Sv_{k}+Sv_{k-1}=\tau ^{2}\psi _{k}\ ,  \label{E1.9}
\end{equation}
where $k=1,\ldots ,n-1$,
\begin{eqnarray*}
S &=&\sum_{j=1}^{m}\eta _{j}S_{j}\ ,\quad \psi _{k}=\sum_{j=1}^{m}\delta
_{1,j}S_{j}\left[ f(t_{k})-M\left( v_{k}\right) \right] \\
&=&S_{1}\left[ f(t_{k})-M\left( v_{k}\right) \right] \ .
\end{eqnarray*}

If we take (\ref{E1.9}) from (\ref{E1.7}), we get
\begin{equation}
z_{k+1}-2Sz_{k}+Sz_{k-1}=r_{k}\ ,  \label{E1.11}
\end{equation}
where $z_{k}=v_{k}-u(t_{k})$,
\begin{eqnarray*}
r_{k} &=&r_{0,k}+r_{1,k}-L\left( \tau ^{2}r_{2,k}+r_{3,k}\right) +r_{4,k}\ ,
\\
r_{0,k} &=&\left( S-L\right) u(t_{k})\ ,\quad r_{1,k}=\left( S-L\right)
\left[ u(t_{k})-u(t_{k-1})\right] \ , \\
r_{2,k} &=&A\left[ u(t_{k+1})-u(t_{k})\right] \ ,\quad r_{4,k}=\tau ^{2}%
\left[ \psi _{k}-L\widetilde{f}(t_{k})\right] \ ,
\end{eqnarray*}

\begin{eqnarray*}
r_{3,k} &=&\int_{t_{k}}^{t_{k+1}}(t_{k+1}-t)\left[ u^{\prime \prime }\left(
t\right) -u^{\prime \prime }\left( t_{k}\right) \right] dt \\
&&+\int_{t_{k-1}}^{t_{k}}(t-t_{k-1})\left[ u^{\prime \prime }\left( t\right)
-u^{\prime \prime }\left( t_{k}\right) \right] dt\ .
\end{eqnarray*}

To present solution of equation (\ref{E1.11}) in explicit form, we
need certain class polynomials, which we call two-variable Chebyshev
polynomials. These polynomials are defined by the following recurrent
relation:
\begin{eqnarray}
\widetilde{U}_{k+1}(x,y) &=&x\widetilde{U}_{k}(x,y)-y\widetilde{U}%
_{k-1}(x,y),\quad k=1,2,...\quad ,  \label{E1.15} \\
\widetilde{U}_{1}(x,y) &=&x,\quad \widetilde{U}_{0}(x,y)=1.  \nonumber
\end{eqnarray}

Notice that the works \cite{R2}, \cite{R3} \ are devoted to research of
three-layer semi-discrete schemes by means of Chebyshev polynomials.

We call $\widetilde{U}_{k}(x,y)$ two-variable Chebyshev polynomials as $%
U_{k}(x)=\widetilde{U}_{k}\left( 2x,1\right) $ represents second order
Chebyshev polynomials.

The following formula can be easily obtained

\begin{equation}
\widetilde{U}_{k}(x,y)=\sqrt{y^{k}}\widetilde{U}_{k}(\xi ,1),\quad \xi =%
\frac{x}{\sqrt{y}},\quad y>0\ ,  \label{E1.21}
\end{equation}

that relates $\widetilde{U}_{k}(x,y)$ to $U_{k}(x)$.

Now we can explicitly write solution of equation (\ref{E1.11}) by means of
polynomials $\widetilde{U}_{k}(x,y)$. Using induction we get
\begin{equation}
z_{k+1}=\widetilde{U}_{k}\left( 2S,S\right) z_{1}-S\widetilde{U}_{k-1}\left(
2S,S\right) z_{0}+\sum_{i=1}^{k}\widetilde{U}_{k-i}\left( 2S,S\right) r_{i}\
.  \label{E1.23}
\end{equation}

Obviously, as $S_{j}$ (j$=1,\ldots ,m$) are self-adjoint, non-negative,
bounded operators, then operator $S$ also will be self-adjoint non-negative
and bounded. Hence, as known (see, e. g., \cite{Reed}, Chapter VI), \ there
exists the only square root $S^{1/2}$ . Taking this into account and using
formula (\ref{E1.21}) , we obtain
\[
\widetilde{U}_{k}\left( 2S,S\right) =B^{k}\widetilde{U}_{k}\left(
2B,I\right) =B^{k}U_{k}(B)\ ,\quad B=S^{1/2}.
\]

Due to this equality, (\ref{E1.23}) will be
\begin{equation}
z_{k+1}=B^{k}U_{k}(B)z_{1}-B^{k+1}U_{k-1}(B)z_{0}+%
\sum_{i=1}^{k}B^{k-i}U_{k-i}(B)r_{i}\ .  \label{E1.25}
\end{equation}

Formula (\ref{E1.25}) is a main relation, by using of which the
following theorem is proved.

\textbf{Theorem 3.1. }If the problem (\ref{2.1}), (\ref{2.1_2}) has a
solution, then the estimate is valid for approximate solution error:

\begin{equation}
\left\Vert z_{k+1}\right\Vert \leq \exp \left( ct_{k-1}\right) \left( \gamma
_{0}\left\Vert \frac{\Delta z_{0}}{\tau }\right\Vert +\gamma _{1}\left\Vert
z_{0}\right\Vert +\Theta _{k}(\tau )\right) \ ,  \label{zk}
\end{equation}
where $z_{k}=v_{k}-u(t_{k})$, $\Delta z_{0}=z_{1}-z_{0}$ , $c=\nu ^{-1/2}a$
, $\nu =\min\limits_{1\leq j\leq m}\left( \alpha _{j}\right) $, $\gamma _{0}=\nu
^{-1/2}+c\tau ^{2}$, $\gamma _{1}=1+c\tau $ ,

\begin{eqnarray*}
\Theta _{k}(\tau ) &=&\tau ^{2}\sum_{i=1}^{k}\left[ c_{1}\left\Vert
Au(t_{i})\right\Vert +c_{3}\left\Vert \widetilde{f}(t_{i})\right\Vert \right]
\\
&&+\tau \sum_{i=1}^{k}\left[ c_{1}J_{i}(t_{i-1},\
A^{1/2}u)+c_{2}J_{i}(t_{i+1},\ A^{1/2}u)\right] \\
&&+c_{2}\sum_{i=1}^{k}\int_{t_{i-1}}^{t_{i+1}}\left[ J_{i}\left( t,\
A^{1/2}u\right) +J_{i}\left( t,\ A^{-1/2}\widetilde{f}\right) \right] dt\ ,
\end{eqnarray*}
and where
\begin{eqnarray*}
J_{i}\left( t,u\right) &=&\left\Vert u(t_{i})-u(t)\right\Vert ,\quad
\widetilde{f}(t)=f(t)-M(u(t)), \\
c_{1} &=&\sum_{j=1}^{m}\eta _{j}^{-3/2}\left( \eta _{j}^{-1}a_{j}+1\right)
,\quad a_{j}=\left\Vert A_{j}A^{-1}\right\Vert <\infty \ , \\
c_{2} &=&m+c_{0}\ ,\quad c_{0}=\sum_{j=1}^{m}\eta _{j}^{-1/2}a_{j}\ ,\quad
c_{3}=\eta _{1}^{-1/2}+m+c_{0}\ .
\end{eqnarray*}

\textbf{Result 2.2.} If the functions $f(t)$ and $A^{1/2}$\bigskip $u(t)$ on
$[0,\ T]$ satisfies Holder condition by the index $\lambda $ ($0<\lambda
\leq 1$), then
\[
\left\Vert u(t_{k})-v_{k}\right\Vert \leq c\tau ^{\lambda },\quad c=const>0.
\]

\setcounter{equation}{0}

\section{Auxiliary lemmas and remarks}

\textbf{Remark 4.1. }As $A_{j}A^{-1}$ ($j=1,\ldots ,m$) are closed operators
(it can be proved easily), therefore due to Closed Graph Theorem they are
bounded, i.e. $a_{j}=\left\Vert A_{j}A^{-1}\right\Vert <\infty $ .

\textbf{Lemma 4.2. }For any $j$ \ ($j=1,\ldots ,m$) $D(A^{1/2})\subset
D(A_{j}^{1/2})$ and
\begin{equation}
\left\Vert A_{j}^{1/2}u\right\Vert \leq \left\Vert A^{1/2}u\right\Vert
,\quad \forall u\in D(A^{1/2}).  \label{E1.27}
\end{equation}

\textbf{Proof: }According to condition (\ref{A.1}) we have:

\[
\left( A_{j}u,u\right) \leq \left( Au,u\right) ,\quad \forall u\in
D(A)\subset D(A_{j}).
\]

From here it follows
\begin{equation}
\left\Vert A_{j}^{1/2}u\right\Vert \leq \left\Vert A^{1/2}u\right\Vert
,\quad \forall u\in D(A)\subset D(A_{j}).  \label{E1.27.1}
\end{equation}

It is known that  $D(A)$ is a ball for $A^{1/2}$(see \cite{TKa}, p. 354). This means: for each $u\in D(A^{1/2})$ there exists a sequence $u_{n}\in D(A)$ such that
 $u_{n}\rightarrow u$ and $A^{1/2}u_{n}\rightarrow A^{1/2}u$.
From here, according to (\ref{E1.27.1}) it follows that $A_{j}^{1/2}u_{n}$
is a Cauchy Sequence. Obviously it will be convergent due to completeness of $H$. Since  $%
A_{j}^{1/2}$ is closed ( $A_{j}$ is given as self-adjoint and positive definite), therefore  $u\in D(A_{j}^{1/2})$ and $%
A_{j}^{1/2}u_{n}\rightarrow A_{j}^{1/2}u$. Thus  $D(A^{1/2})\subset
D(A_{j}^{1/2})$, and in addition, the inequality
\[
\left\Vert A_{j}^{1/2}u_{n}\right\Vert \leq \left\Vert
A^{1/2}u_{n}\right\Vert ,\quad u_{n}\in D(A)\subset D(A_{j})
\]%
gives (\ref{E1.27}).

\textbf{Remark 4.3. }If self-adjoint positive definite operators $A$ and $B$
are such that $D(A)\subset D(B)$ and $B\leq A$ ($\left( Bu,u\right) \leq
\left( Au,u\right) $ , $\forall u\in D(A)$ ), then $A^{-1}\leq B^{-1}$ .

Indeed, let's  $u=B^{-1}f$ \ and $\ v=A^{-1}f$ , $f\in H$. The relation is true (see the proof of Theorem VI.2.21 from \cite{TKa}):
\begin{eqnarray*}
\left( A^{-1}f,\ f\right) ^{2} &=&\left( v,\ Bu\right) ^{2}=\left(
B^{1/2}v,\ B^{1/2}u\right) ^{2} \\
&\leq &\left\Vert B^{1/2}v\right\Vert ^{2}\left\Vert B^{1/2}u\right\Vert
^{2}=\left( Bv,v\right) \left( Bu,u\right) \\
&\leq &\left( Av,v\right) \left( Bu,u\right) =\left( f,\ A^{-1}f\right)
\left( f,\ B^{-1}f\right) \\
&=&\left( A^{-1}f,\ f\right) \left( B^{-1}f,\ f\right) .
\end{eqnarray*}

After reduction we get
\[
\left( A^{-1}f,\ f\right) \leq \left( B^{-1}f,\ f\right) ,
\]
i. e. $A^{-1}\leq B^{-1}$.

\textbf{Remark 4.4 }The formula is valid
\begin{equation}
S-L=\tau ^{2}\sum_{j=1}^{m}\eta _{j}^{-1}\left( I-S_{j}\right) \left( \eta
_{j}^{-1}A_{j}A^{-1}-I\right) AL\ .  \label{E1.35}
\end{equation}

Indeed, as
\[
I-S=\sum_{j=1}^{m}\eta _{j}\left( I-S_{j}\right) =\tau
^{2}\sum_{j=1}^{m}A_{j}S_{j}\ ,
\]
therefore
\begin{eqnarray*}
S-L &=&\left[ S\left( I+\tau ^{2}A\right) -I\right] L=\left[ \left(
S-I\right) +\tau ^{2}AS\right] L \\
&=&\tau ^{2}\sum_{j=1}^{m}\left( \eta _{j}S_{j}A-A_{j}S_{j}\right) L\ .
\end{eqnarray*}

From here, taking the equality into account
\begin{eqnarray*}
\eta _{j}S_{j}A-A_{j}S_{j} &=&\eta _{j}S_{j}A-S_{j}A_{j}=S_{j}\left( \eta
_{j}A-A_{j}\right) \\
&=&\left( S_{j}-I\right) \left( \eta _{j}A-A_{j}\right) +\left( \eta
_{j}A-A_{j}\right) \ ,
\end{eqnarray*}
we get (\ref{E1.35}) .

We can consider the closeness of operators $S$ and $L$ also by the
following formulas:
\begin{eqnarray*}
S_{j} &=&\tau ^{4}\eta _{j}^{-2}A_{j}^{2}S_{j}-\tau ^{2}\eta
_{j}^{-1}A_{j}+I\ , \\
S &=&\tau ^{4}\sum_{j=1}^{m}\eta _{j}^{-1}A_{j}^{2}S_{j}-\tau ^{2}A+I\ , \\
L &=&\tau ^{4}A^{2}L-\tau ^{2}A+I\ .
\end{eqnarray*}

\textbf{Lemma 4.5. }The inequalities are valid:
\begin{eqnarray}
&&\left\Vert \left( I-S\right) ^{-1/2}ALf\right\Vert  \nonumber \\
&\leq &\tau ^{-1}\left( m\left\Vert A^{1/2}Lf\right\Vert +c_{0}\left\Vert
A^{1/2}L^{1/2}f\right\Vert \right) ,\quad f\in H,  \label{E1.39}
\end{eqnarray}

\begin{equation}
\left\Vert \left( I-S\right) ^{-1/2}\left( S-L\right) f\right\Vert \leq \tau
^{2}c_{1}\left\Vert ALf\right\Vert ,\quad f\in H,  \label{E1.43}
\end{equation}

\begin{equation}
\left\Vert \left( I-S\right) ^{-1/2}LAu\right\Vert \leq \tau
^{-1}c_{2}\left\Vert A^{1/2}u\right\Vert ,\quad u\in D(A),  \label{E1.45}
\end{equation}

\begin{equation}
\left\Vert \left( I-S\right) ^{-1/2}Lf\right\Vert \leq \tau
^{-1}c_{2}\left\Vert A^{-1/2}f\right\Vert ,\quad f\in H.  \label{E1.47}
\end{equation}

\textbf{Proof: }Let us prove the inequality (\ref{E1.39}) . As
\[
I-S\geq \eta _{j}\left( I-S_{j}\right) =\tau ^{2}A_{j}S_{j}>0,
\]
therefore (see remark 4.3)
\begin{eqnarray}
\left( I-S\right) ^{-1} &\leq &\eta _{j}^{-1}\left( I-S_{j}\right)
^{-1}=\tau ^{-2}\left( I+\tau ^{2}\eta _{j}^{-1}A_{j}\right) A_{j}^{-1}
\nonumber \\
&=&\tau ^{-2}A_{j}^{-1}+\eta _{j}^{-1}I\leq \left( \tau ^{-1}A_{j}^{-s}+\eta
_{j}^{-s}I\right) ^{2},\quad s=\frac{1}{2}\ .  \label{E1.48}
\end{eqnarray}

Thus we have
\begin{eqnarray*}
\left\Vert \left( I-S\right) ^{-s}f\right\Vert &\leq &\left\Vert \left( \tau
^{-1}A_{j}^{-s}+\eta _{j}^{-s}I\right) f\right\Vert \\
&\leq &\tau ^{-1}\left\Vert A_{j}^{-s}f\right\Vert +\eta _{j}^{-s}\left\Vert
f\right\Vert ,\quad f\in H.
\end{eqnarray*}

Using this inequality we get (below everywhere $s=1/2$):
\begin{eqnarray}
\left\Vert \left( I-S\right) ^{-s}ALf\right\Vert &=&\left\Vert \left(
I-S\right) ^{-s}\sum_{j=1}^{m}A_{j}Lf\right\Vert  \nonumber \\
&\leq &\tau ^{-1}\sum_{j=1}^{m}\left( \left\Vert A_{j}^{s}Lf\right\Vert
+\tau \eta _{j}^{-s}\left\Vert A_{j}Lf\right\Vert \right) .  \label{E1.49}
\end{eqnarray}

According to lemma 4.2 we have
\begin{equation}
\left\Vert A_{j}^{s}Lf\right\Vert \leq \left\Vert A^{s}Lf\right\Vert .
\label{E1.51}
\end{equation}

According to remark 4.1 we have
\begin{equation}
\left\Vert A_{j}Lf\right\Vert \leq \left\Vert A_{j}A^{-1}\right\Vert
\left\Vert ALf\right\Vert =a_{j}\left\Vert ALf\right\Vert .  \label{E1.53}
\end{equation}

From (\ref{E1.49}), taking into account (\ref{E1.51}) and (\ref{E1.53}), we
obtain
\[
\left\Vert \left( I-S\right) ^{-s}ALf\right\Vert \leq m\tau ^{-1}\left\Vert
A^{s}Lf\right\Vert +c_{0}\left\Vert ALf\right\Vert .
\]

Hence, taking into account inequality
\[
\left\Vert ALf\right\Vert =\left\Vert A^{s}L^{s}\left( A^{s}L^{s}f\right)
\right\Vert \leq \left\Vert A^{s}L^{s}\right\Vert \left\Vert
A^{s}L^{s}f\right\Vert \leq \tau ^{-1}\left\Vert A^{s}L^{s}f\right\Vert
\]

(\ref{E1.39}) is obtained.

Let us prove inequality (\ref{E1.43}). From (\ref{E1.48}) we have

\begin{equation}
\left\Vert \left( I-S\right) ^{-1/2}h\right\Vert \leq \eta
_{j}^{-1/2}\left\Vert \left( I-S_{j}\right) ^{-1/2}h\right\Vert ,\quad h\in
H.  \label{E1.55}
\end{equation}

If in (\ref{E1.55}) we substitute $h=\left( I-S_{j}\right) f$ , we get
\begin{equation}
\left\Vert \left( I-S\right) ^{-1/2}\left( I-S_{j}\right) f\right\Vert \leq
\eta _{j}^{-1/2}\left\Vert \left( I-S_{j}\right) ^{1/2}f\right\Vert \leq
\eta _{j}^{-1/2}\left\Vert f\right\Vert .  \label{E1.59}
\end{equation}

From this inequality and (\ref{E1.35}) it follows
\begin{eqnarray*}
&&\left\Vert \left( I-S\right) ^{-1/2}\left( S-L\right) f\right\Vert \\
&\leq &\tau ^{2}\sum_{j=1}^{m}\eta _{j}^{-1}\left\Vert \left( I-S\right)
^{-1/2}\left( I-S_{j}\right) \left( \eta _{j}^{-1}A_{j}A^{-1}-I\right)
ALf\right\Vert \\
&\leq &\tau ^{2}\sum_{j=1}^{m}\eta _{j}^{-3/2}\left\Vert \left( \eta
_{j}^{-1}A_{j}A^{-1}-I\right) ALf\right\Vert \leq \tau ^{2}c_{1}\left\Vert
ALf\right\Vert \ .
\end{eqnarray*}

Let us prove the inequality (\ref{E1.45}) . Obviously, operators $A$ and $L$
are commutative on $D(A)$. Then due to (\ref{E1.39}), the relation is valid:
\begin{eqnarray*}
\left\Vert \left( I-S\right) ^{-1/2}LAu\right\Vert &=&\left\Vert \left(
I-S\right) ^{-1/2}ALu\right\Vert \\
&\leq &\tau ^{-1}\left( m\left\Vert A^{1/2}Lu\right\Vert +c_{0}\left\Vert
A^{1/2}L^{1/2}u\right\Vert \right) \\
&\leq &\tau ^{-1}\left( m\left\Vert LA^{1/2}u\right\Vert +c_{0}\left\Vert
L^{1/2}A^{1/2}u\right\Vert \right) \\
&\leq &\tau ^{-1}c_{2}\left\Vert A^{1/2}u\right\Vert .
\end{eqnarray*}

Obviously inequality (\ref{E1.47}) is a result of (\ref{E1.45}).

\textbf{Lemma 4.6. }For operator polynomials $U_{k}(B)$ the following
estimate is valid:
\begin{equation}
\left\Vert BU_{k}(B)\right\Vert \leq \left( \tau \sqrt{\nu }\right)
^{-1},\quad \nu =\min_{1\leq j\leq m}\left( \alpha _{j}\right) \ ,
\label{B1}
\end{equation}

\begin{equation}
\left\Vert U_{k}(B)\left( I-B^{2}\right) ^{1/2}\right\Vert \leq 1\ ,
\label{B2}
\end{equation}

\begin{equation}
\left\Vert U_{k}(B)-BU_{k-1}(B)\right\Vert \leq 1,  \label{B3}
\end{equation}

\begin{equation}
\left\Vert BU_{k}(B)-U_{k-1}(B)\right\Vert \leq 1 .  \label{B4}
\end{equation}

\textbf{Proof: }As it is well-known, for Chebyshev second-order polynomials the
estimate is valid (see, e. g., \cite{Sege}):

\begin{equation}
\left\vert U_{k}(x)\right\vert \leq \frac{1}{\sqrt{1-x^{2}}}\ ,\quad x\in
]-1,1[\ ,  \label{U1}
\end{equation}

that follows from the well-known formula
\[
U_{k}\left( x\right) =\frac{\sin \left( \left( k+1\right) \arccos x\right) }{%
\sqrt{1-x^{2}}},\quad x\in ]-1,1[\ .
\]

From this formula, by means of simple calculations, we also obtain the
estimate
\begin{equation}
\left\vert U_{k}(x)-U_{k-1}(x)\right\vert \leq \sqrt{\frac{2}{1+x}}\ ,\quad
x\in ]-1,1[\ .  \label{U3}
\end{equation}

Indeed we have
\begin{eqnarray*}
\left\vert U_{k}(x)-U_{k-1}(x)\right\vert &=&\frac{2}{\sqrt{1-x^{2}}}%
\left\vert \cos (k+\frac{\theta }{2})\sin \frac{\theta }{2}\right\vert \\
&\leq &\frac{2}{\sqrt{1-x^{2}}}\sqrt{\frac{1-\cos \theta }{2}} \\
&=&\sqrt{\frac{2}{1+x}}\ ,\quad x\in ]-1,1[\ ,
\end{eqnarray*}
where $\theta =\arccos x$ .

The following estimation is valid:
\begin{equation}
\left\vert U_{k}(x)-xU_{k-1}(x)\right\vert  \leq  \left\vert \cos((k+1)\theta)\right\vert \leq 1.  \label{U5}
\end{equation}

Obviously, similarly we get the estimate
\begin{equation}
\left\vert xU_{k}(x)-U_{k-1}(x)\right\vert \leq 1 .  \label{U7}
\end{equation}

For the estimation of norm of polynomial operators $U_{k}(B)$ , we
need to estimate the spectrum of operator $B=S^{1/2}$, that obviously can be
reduced to the estimate of spectrum of operator $S$ .

Let us first estimate spectrum of operator $S_{j}=\left( I+\tau
^{2}\eta _{j}^{-1}A_{j}\right) ^{-1}$ . We obviously have
\[
\left( I+\tau ^{2}\eta _{j}^{-1}A_{j}\right) \geq \left( 1+\tau ^{2}\eta
_{j}^{-1}\alpha _{j}\right) I>0\ .
\]

From here, due to remark 4.3, it follows
\begin{equation}
0<S_{j}\leq \left( 1+\tau ^{2}\eta _{j}^{-1}\alpha _{j}\right) ^{-1}I\leq
\left( 1+\tau ^{2}\nu \right) ^{-1}I\ .  \label{U9}
\end{equation}

If we take into account representation of $S$, then according to (\ref{U9}),
we have
\[
0<S\leq \sum_{j=1}^{m}\eta _{j}\left( 1+\tau ^{2}\nu \right) ^{-1}I=\left(
1+\tau ^{2}\nu \right) ^{-1}I\ .
\]

This means that
\[
Sp\left( S\right) \subset \lbrack 0,\ \left( 1+\tau ^{2}\nu \right) ^{-1}]
\]

From here it follows that (due to well-known theorem on spectral mapping)
\begin{equation}
Sp\left( B\right) \subset \lbrack 0,\ \left( 1+\tau ^{2}\nu \right)
^{-1/2}]\ .  \label{SpB}
\end{equation}

Further we can easily show estimates of (\ref{B1})-(\ref{B4}).

Let us show estimate (\ref{B1}). As is known, norm of
operator-function, when the argument represents self-adjoint bounded
operator, is equal to the $C$-norm of the corresponding scalar function
on the spectrum (see, e. g., \cite{Reed}, Chapter VII). Due to this result
we have
\[
\left\Vert BU_{k}(B)\right\Vert \leq \max_{x\in Sp(B)}\left\vert
xU_{k}(x)\right\vert \ .
\]

From here, by the inequality (\ref{U1}) and relation (\ref{SpB}), we get
\[
\left\Vert BU_{k}(B)\right\Vert \leq \max_{x\in Sp(B)}\frac{x}{\sqrt{1-x^{2}}%
}\leq \frac{1}{\tau \sqrt{\nu }}\ .
\]

Similarly is obtained (\ref{B2}). Indeed we have
\begin{eqnarray*}
\left\Vert U_{k}(B)\left( I-B^{2}\right) ^{1/2}\right\Vert &\leq &\max_{x\in
Sp(B)}\left\vert U_{k}(x)\sqrt{1-x^{2}}\right\vert \\
&\leq &\max_{x\in Sp(B)}\left( \frac{1}{\sqrt{1-x^{2}}}\sqrt{1-x^{2}}\right)
=1\ .
\end{eqnarray*}

From the inequalities (\ref{U5}) and (\ref{U7}) respectively follows the
estimates (\ref{B3}) and (\ref{B4}) .

\setcounter{equation}{0}

\section{Error estimate of approximate solution}

In this section we continue to prove Theorem 3.1, that concerns error
estimate of approximate solution.

Let us rewrite formula (\ref{E1.25}) in the following form%
\begin{eqnarray}
z_{k+1} &=&\tau B^{k}U_{k}(B)\frac{\Delta z_{0}}{\tau }+B^{k}\left(
U_{k}(B)-BU_{k-1}(B)\right) z_{0}  \nonumber \\
&&+\sum_{i=1}^{k}B^{k-i}U_{k-i}(B)r_{i}\ ,  \label{zk1.1}
\end{eqnarray}
where $\Delta z_{0}=z_{1}-z_{0}$ .

Let us note that $\Delta z_{0}/\tau $ is an error of the difference
analog of the first order derivative of approximate solution at the point $%
t=0 $ ,

\[
\frac{\Delta z_{0}}{\tau }=\frac{\Delta u_{0}}{\tau }-\frac{\Delta u(0)}{%
\tau }\ ,\quad \Delta u(0)=u(\tau )-u(0)\ .
\]

If we move to norms in (\ref{zk1.1}), we obtain
\begin{eqnarray*}
\left\Vert z_{k+1}\right\Vert &\leq &\tau \left\Vert B^{k-1}\right\Vert
\left\Vert BU_{k}(B)\right\Vert \left\Vert \frac{\Delta z_{0}}{\tau }%
\right\Vert \\
&&+\left\Vert B^{k}\right\Vert \left\Vert U_{k}(B)-BU_{k-1}(B)\right\Vert
\left\Vert z_{0}\right\Vert \\
&&+\sum_{i=1}^{k}\left\Vert B^{k-i}\right\Vert \left\Vert U_{k-i}(B)\left(
I-B^{2}\right) ^{1/2}\right\Vert \left\Vert \left( I-B^{2}\right)
^{-1/2}r_{i}\right\Vert \ .
\end{eqnarray*}

From here, taking into account estimates (\ref{B1}), (\ref{B2}) and (\ref{B3}%
) ($\left\Vert B\right\Vert \leq 1$), we have
\begin{equation}
\left\Vert z_{k+1}\right\Vert \leq \nu ^{-1/2}\left\Vert \frac{\Delta z_{0}}{%
\tau }\right\Vert +\left\Vert z_{0}\right\Vert
+\sum_{i=1}^{k}\left\Vert \left( I-B^{2}\right) ^{-1/2}r_{i}\right\Vert \ .
\label{zk1.5}
\end{equation}

Obviously we have
\begin{eqnarray}
\left\Vert \left( I-B^{2}\right) ^{-1/2}r_{i}\right\Vert &=&\left\Vert
\left( I-S\right) ^{-1/2}r_{i}\right\Vert  \nonumber \\
&\leq &\lambda _{0,i}+\lambda _{1,i}+\tau ^{2}\lambda _{2,i}+\lambda
_{3,i}+\lambda _{4,i}\ ,  \label{RI}
\end{eqnarray}
where
\begin{eqnarray*}
\lambda _{s,i} &=&\left\Vert \left( I-S\right) ^{-1/2}r_{s,i}\right\Vert \
,\quad s=0,\ 1,\ 4, \\
\lambda _{s,i} &=&\left\Vert \left( I-S\right) ^{-1/2}Lr_{s,i}\right\Vert \
,\quad s=2,\ 3\ .
\end{eqnarray*}

Using (\ref{E1.43}), we get:
\begin{eqnarray}
\lambda _{0,i} &=&\left\Vert \left( I-S\right) ^{-1/2}\left( S-L\right)
u(t_{i})\right\Vert  \nonumber \\
&\leq &\tau ^{2}c_{1}\left\Vert ALu(t_{i})\right\Vert \leq \tau
^{2}c_{1}\left\Vert Au(t_{i})\right\Vert ,  \label{lamda0}
\end{eqnarray}

\begin{eqnarray}
\lambda _{1,i} &=&\left\Vert \left( I-S\right) ^{-1/2}\left( S-L\right)
\left[ u(t_{i})-u(t_{i-1})\right] \right\Vert  \nonumber \\
&\leq &\tau ^{2}c_{1}\left\Vert AL\left[ u(t_{i})-u(t_{i-1})\right]
\right\Vert  \nonumber \\
&\leq &\tau c_{1}\left\Vert \tau A^{1/2}L\right\Vert \left\Vert A^{1/2}\left[
u(t_{i})-u(t_{i-1})\right] \right\Vert  \nonumber \\
&\leq &\tau c_{1}\left\Vert A^{1/2}\left[ u(t_{i})-u(t_{i-1})\right]
\right\Vert .  \label{lamda1}
\end{eqnarray}

Using (\ref{E1.45}) and (\ref{E1.47}), respectively we get:
\begin{eqnarray}
\tau ^{2}\lambda _{2,i} &=&\tau ^{2}\left\Vert \left( I-S\right) ^{-1/2}LA
\left[ u(t_{i+1})-u(t_{i})\right] \right\Vert  \nonumber \\
&\leq &\tau c_{2}\left\Vert A^{1/2}\left[ u(t_{i+1})-u(t_{i})\right]
\right\Vert ,  \label{lamda2}
\end{eqnarray}

\[
\lambda _{3,i}=\left\Vert \left( I-S\right) ^{-1/2}Lr_{3,i}\right\Vert \leq
\tau ^{-1}c_{2}\left\Vert A^{-1/2}r_{3,i}\right\Vert
\]

\begin{eqnarray*}
&\leq &\tau ^{-1}c_{2}\left\Vert \int_{t_{i}}^{t_{i+1}}(t_{i+1}-t)A^{-1/2}
\left[ u^{\prime \prime }\left( t\right) -u^{\prime \prime }\left(
t_{i}\right) \right] dt\right\Vert \\
&&+\tau ^{-1}c_{2}\left\Vert \int_{t_{i-1}}^{t_{i}}(t-t_{i-1})A^{-1/2}\left[
u^{\prime \prime }\left( t\right) -u^{\prime \prime }\left( t_{i}\right) %
\right] dt\right\Vert
\end{eqnarray*}

\[
\leq c_{2}\int_{t_{i-1}}^{t_{i+1}}\left\Vert A^{-1/2}\left[ u^{\prime \prime
}\left( t\right) -u^{\prime \prime }\left( t_{i}\right) \right] \right\Vert
dt\ .
\]

From here, taking into account equation (\ref{2.1}) , we have
\begin{eqnarray}
\lambda _{3,i} &\leq &c_{2}\int_{t_{i-1}}^{t_{i+1}}\left\Vert A^{1/2}\left[
u\left( t\right) -u\left( t_{i}\right) \right] \right\Vert dt  \nonumber \\
&&+c_{2}\int_{t_{i-1}}^{t_{i+1}}\left\Vert A^{-1/2}\left[ \widetilde{f}(t)-%
\widetilde{f}(t_{i})\right] \right\Vert dt\ .  \label{lamda3}
\end{eqnarray}

Let us estimate $\lambda _{4,i}$ . Let's rewrite $r_{4,i}$ as
\begin{eqnarray}
r_{4,i} &=&\tau ^{2}\left[ \psi _{i}-L\widetilde{f}(t_{i})\right] =\tau
^{2}S_{1}\left[ f(t_{i})-M\left( v_{i}\right) \right] -\tau ^{2}L\widetilde{f%
}(t_{i})  \nonumber \\
&=&\tau ^{2}S_{1}\left[ f(t_{i})-M(u(t_{i}))\right] +\tau ^{2}S_{1}\left[
M(u(t_{i}))-M\left( v_{i}\right) \right] -\tau ^{2}L\widetilde{f}(t_{i})
\nonumber \\
&=&\tau ^{2}S_{1}\left[ M(u(t_{i}))-M\left( v_{i}\right) \right] +\tau
^{2}S_{1}\widetilde{f}(t_{i})-\tau ^{2}L\widetilde{f}(t_{i})  \nonumber \\
&=&\tau ^{2}S_{1}\left[ M(u(t_{i}))-M\left( v_{i}\right) \right] +\tau
^{2}(S_{1}-I)\widetilde{f}(t_{i})+\tau ^{2}(I-L)\widetilde{f}(t_{i})
\nonumber \\
&=&\tau ^{2}S_{1}\left[ M(u(t_{i}))-M\left( v_{i}\right) \right] +\tau
^{2}(S_{1}-I)\widetilde{f}(t_{i})+\tau ^{4}AL\widetilde{f}(t_{i})\ .
\label{R4}
\end{eqnarray}

Hence we have
\begin{eqnarray}
\lambda _{4,i} &=&\left\Vert \left( I-S\right) ^{-1/2}r_{4,i}\right\Vert
\nonumber \\
&\leq &\tau ^{2}\left\Vert \left( I-S\right) ^{-1/2}S_{1}\left[
M(u(t_{i}))-M\left( v_{i}\right) \right] \right\Vert  \nonumber \\
&&+\tau ^{2}\left\Vert \left( I-S\right) ^{-1/2}\left( I-S_{1}\right)
\widetilde{f}(t_{i})\right\Vert  \nonumber \\
&&+\tau ^{4}\left\Vert \left( I-S\right) ^{-1/2}AL\widetilde{f}%
(t_{i})\right\Vert .  \label{lamda4.1}
\end{eqnarray}

If we use inequality (\ref{E1.55}), we get:
\begin{eqnarray}
\left\Vert \left( I-S\right) ^{-1/2}S_{1}h\right\Vert &\leq &\eta
_{1}^{-1/2}\left\Vert \left( I-S_{1}\right) ^{-1/2}S_{1}h\right\Vert
\nonumber \\
&=&\eta _{1}^{-1/2}\left\Vert \left( \tau ^{2}\eta
_{1}^{-1}A_{1}S_{1}\right) ^{-1/2}S_{1}h\right\Vert  \nonumber \\
&=&\tau ^{-1}\left\Vert A_{1}^{-1/2}S_{1}^{1/2}h\right\Vert \leq \left( \tau
\sqrt{\nu }\right) ^{-1}\left\Vert h\right\Vert ,  \label{lamda4.2}
\end{eqnarray}

\begin{equation}
\left\Vert \left( I-S\right) ^{-1/2}\left( I-S_{1}\right) h\right\Vert \leq
\eta _{1}^{-1/2}\left\Vert \left( I-S_{1}\right) ^{1/2}h\right\Vert \leq
\eta _{1}^{-1/2}\left\Vert h\right\Vert .  \label{lamda4.3}
\end{equation}

Taking into account that$\left\Vert \tau A^{1/2}L^{1/2}\right\Vert \leq 1$,
then from (\ref{E1.39}) we get
\begin{eqnarray}
&&\tau ^{2}\left\Vert \left( I-S\right) ^{-1/2}ALf\right\Vert  \nonumber \\
&\leq &\tau \left( m\left\Vert A^{1/2}Lf\right\Vert +c_{0}\left\Vert
A^{1/2}L^{1/2}f\right\Vert \right) \leq (m+c_{0})\left\Vert f\right\Vert .
\label{lamda4.4}
\end{eqnarray}

From (\ref{lamda4.1}), taking into account estimates (\ref{lamda4.2}), (\ref%
{lamda4.3}) and (\ref{lamda4.4}) , we get
\begin{equation}
\lambda _{4,i}\leq \tau \nu ^{-1/2}\left\Vert M(u(t_{i}))-M\left(
v_{i}\right) \right\Vert +\tau ^{2}c_{3}\left\Vert \widetilde{f}%
(t_{i})\right\Vert .  \label{lamda4.5}
\end{equation}
where $c_{3}=\eta _{1}^{-1/2}+m+c_{0}$ .

Since nonlinear operator $M$ satisfies Lipschietz condition,
therefore from (\ref{lamda4.5}) we have

\begin{equation}
\lambda _{4,i}\leq \tau \nu ^{-1/2}a\left\Vert u(t_{i})-v_{i}\right\Vert
+\tau ^{2}c_{3}\left\Vert \widetilde{f}(t_{i})\right\Vert ,  \label{lamda4.6}
\end{equation}
where $a$ is a Lipchietz constant.

If in (\ref{RI}) we substitute (\ref{lamda0}), (\ref{lamda1}), (\ref%
{lamda2}), (\ref{lamda3}) and (\ref{lamda4.6}), we get
\begin{eqnarray}
&&\left\Vert \left( I-B^{2}\right) ^{-1/2}r_{i}\right\Vert  \nonumber \\
&\leq &\tau \nu ^{-1/2}a\left\Vert z_{i}\right\Vert +\tau
^{2}c_{1}\left\Vert Au(t_{i})\right\Vert +\tau c_{1}J_{i}(t_{i-1},\ A^{1/2}u)
\nonumber \\
&&+\tau c_{2}J_{i}(t_{i+1},\ A^{1/2}u)+\tau ^{2}c_{3}\left\Vert \widetilde{f}%
(t_{i})\right\Vert  \nonumber \\
&&+c_{2}\int_{t_{i-1}}^{t_{i+1}}\left[ J_{i}\left( t,\ A^{1/2}u\right)
+J_{i}\left( t,\ A^{-1/2}\widetilde{f}\right) \right] dt\ ,  \label{RI1}
\end{eqnarray}
where
\[
J_{i}\left( t,u\right) =\left\Vert u(t_{i})-u(t)\right\Vert .
\]

From (\ref{zk1.5}), taking into account (\ref{RI1}),we get
\begin{equation}
\delta _{k+1}\leq c\tau \delta _{k}+\lambda _{k}\ ,  \label{deltak}
\end{equation}
where $\delta _{k}=\left\Vert z_{k}\right\Vert $ , $c=\nu ^{-1/2}a$ ,
\[
\lambda _{k}=\nu ^{-1/2}\left\Vert \frac{\Delta z_{0}}{\tau }\right\Vert
+ \left\Vert z_{0}\right\Vert +\Theta _{k}(\tau )\ .
\]

From (\ref{deltak}), according to discrete analog of Gronwell's lemma, we
have
\[
\delta _{k+1}\leq \exp \left( ct_{k-1}\right) \left( c\tau \delta
_{1}+\lambda _{k}\right) \ .
\]

Obviously from here, taking into account inequality
\[
\delta _{1}=\left\Vert z_{1}\right\Vert \leq \tau \left\Vert \frac{\Delta
z_{0}}{\tau }\right\Vert +\left\Vert z_{0}\right\Vert
\]

We get the estimate (\ref{zk}) .

\setcounter{equation}{0}

\section{ Error estimate for the difference analog of the first order
derivative of approximate solution}

In this section, on the basis of the results obtained in the previous
section, we obtain the a priori estimates for error of the difference analog
of the first order derivative of the decomposition (\ref{E1.1}) scheme
solution.

\textbf{Theorem 6.1. }If the problem (\ref{2.1}), (\ref{2.1_2}) has a
solution and $\varphi _{0}\in D(A)$, then the estimate is valid:

\begin{equation}
\left\Vert \frac{\Delta z_{k}}{\tau }\right\Vert \leq \left\Vert
\frac{\Delta z_{0}}{\tau }\right\Vert +\left( \tau t_{k}\right)
^{-1/2}\left\Vert z_{0}\right\Vert +\widetilde{\Theta }_{k}(\tau )\ ,
\label{zkdel}
\end{equation}

where $\Delta z_{k}=z_{k+1}-z_{k}$.

\begin{eqnarray*}
\widetilde{\Theta }_{k}(\tau ) &=&\tau \sum_{i=1}^{k}\left[ c_{5}J_{i}\left(
t_{i-1},u^{\prime \prime }\right) +J_{i}\left( t_{i+1},u^{\prime
\prime }\right) +ac_{6}J_{i}\left( t_{i-1},u\right) \right] \\
&&+\tau \sum_{i=1}^{k}\left[ c_{7}J_{i}\left( t_{i-1},f\right) + J_{i}\left( t_{i-1},f\right) \right] \\
&&+\sum_{i=1}^{k}\left( \int_{t_{i-1}}^{t_{i+1}}J_{i}\left(
t,u^{\prime \prime }\right) dt+\tau a\left\Vert z_{i}\right\Vert \right) \\
&&+\tau \left[ c_{1}\left\Vert A\varphi _{0}\right\Vert +c_{6}\left(
\left\Vert f(0)\right\Vert +\left\Vert M(\varphi _{0})\right\Vert \right) %
\right] \ ,
\end{eqnarray*}
and where $c_{5}=c_{4}+c_{1}$ , $\ c_{6}=\eta _{1}^{-1/2}+m+c_{0}$ ,
$c_{7}=c_{4}+c_{1}+c_{6}$ ,

\[
c_{4}=\sum_{j=1}^{m}\eta _{j}^{-1}\left( \eta _{j}^{-1}a_{j}+1\right) \
,\quad J_{i}\left( t,u\right) =\left\Vert u(t_{i})-u(t)\right\Vert \ .
\]

\textbf{Proof: }Due to formula (\ref{E1.25}) we have
\begin{equation}
z_{k+1}-z_{k}=\left( R_{k}-R_{k-1}\right) z_{1}-\left(
R_{k-1}-R_{k-2}\right) B^{2}z_{0}+\Phi _{k}\ ,  \label{zk2.1}
\end{equation}
where $R_{k}=B^{k}U_{k}(B)$ ,
\[
\Phi _{k}=\sum_{i=1}^{k}R_{k-i}r_{i}-\sum_{i=1}^{k-1}R_{k-1-i}r_{i}\ .
\]

Let us rewrite (\ref{zk2.1}) as
\begin{eqnarray}
\frac{\Delta z_{k}}{\tau } &=&\left( R_{k}-R_{k-1}\right) \frac{\Delta z_{0}%
}{\tau }  \nonumber \\
&&+\tau ^{-1}\left[ R_{k}+B^{2}R_{k-2}-\left( I+B^{2}\right) R_{k-1}\right]
z_{0}+\tau ^{-1}\Phi _{k}\ ,  \label{zk2.3}
\end{eqnarray}
where $\Delta z_{k}=z_{k+1}-z_{k}$ .

By the simple transformation we get
\begin{eqnarray*}
&&R_{k}+B^{2}R_{k-2}-\left( I+B^{2}\right) R_{k-1} \\
&=&B^{k-1}\left[ B\left( U_{k}+U_{k-2}\right) -\left( I+B^{2}\right) U_{k-1}%
\right] \\
&=&-B^{k-1}\left( I-B^{2}\right) U_{k-1}\ .\
\end{eqnarray*}

From here, taking into account (\ref{B2}), we have
\begin{eqnarray}
&&\left\Vert R_{k}+B^{2}R_{k-2}-\left( I+B^{2}\right) R_{k-1}\right\Vert
\nonumber \\
&\leq &\left\Vert B^{k-1}\left( I-B^{2}\right) ^{1/2}\right\Vert \left\Vert
U_{k-1}\left( I-B^{2}\right) ^{1/2}\right\Vert \leq \left\Vert B^{k-1}\left(
I-B^{2}\right) ^{1/2}\right\Vert  \nonumber \\
&\leq &\max_{0\leq x\leq 1}\left[ x^{k-1}\left( 1-x^{2}\right) ^{1/2}\right]
\leq \frac{1}{\sqrt{k}}\ .  \label{zk2.5}
\end{eqnarray}

Obviously, according to (\ref{B4}), we have
\begin{eqnarray}
\left\Vert R_{k}-R_{k-1}\right\Vert &=&\left\Vert B^{k-1}\left(
BU_{k}(B)-U_{k-1}(B)\right) \right\Vert  \nonumber \\
&\leq &\left\Vert B^{k-1}\right\Vert \left\Vert
BU_{k}(B)-U_{k-1}(B)\right\Vert \leq 1 .  \label{zk2.7}
\end{eqnarray}

We can give to $\Phi _{k}$ the following form
\begin{eqnarray*}
\Phi _{k} &=&\sum_{i=1}^{k}R_{k-i}\left( r_{i}-r_{0,i}-\tau ^{2}\zeta
_{i}\right) -\sum_{i=1}^{k-1}R_{k-1-i}\left( r_{i}-r_{0,i}-\tau ^{2}\zeta
_{i}\right) \\
&&+\sum_{i=1}^{k}R_{k-i}\left( r_{0,i}+\tau ^{2}\zeta _{i}\right)
-\sum_{i=1}^{k-1}R_{k-1-i}\left( r_{0,i}+\tau ^{2}\zeta _{i}\right)
\end{eqnarray*}

\begin{eqnarray*}
&=&\sum_{i=1}^{k}\left( R_{k-i}-R_{k-1-i}\right) \left( r_{i}-r_{0,i}-\tau
^{2}\zeta _{i}\right) \\
&&+\sum_{i=1}^{k}R_{k-i}\left( r_{0,i}+\tau ^{2}\zeta _{i}\right)
-\sum_{i=2}^{k}R_{k-i}\left( r_{0,i-1}+\tau ^{2}\zeta _{i-1}\right)
\end{eqnarray*}

\begin{eqnarray}
&=&\sum_{i=1}^{k}\left( R_{k-i}-R_{k-1-i}\right) \left( r_{i}-r_{0,i}-\tau
^{2}\zeta _{i}\right)  \nonumber \\
&&+\sum_{i=1}^{k}R_{k-i}\left[ \left( r_{0,i}-r_{0,i-1}\right) +\tau
^{2}\left( \zeta _{i}-\zeta _{i-1}\right) \right]  \nonumber \\
&&+R_{k-1}\left( r_{0,0}+\tau ^{2}\zeta _{0}\right) \ ,  \label{zk2.15}
\end{eqnarray}
where $R_{-1}=0$ ,
\[
\zeta _{i}=(S_{1}-I)\widetilde{f}(t_{i})+\tau ^{2}AL\widetilde{f}(t_{i})\ .
\]

Taking into account (\ref{zk2.7}), we have
\begin{eqnarray}
&&\left\Vert \left( R_{k-i}-R_{k-1-i}\right) \left( r_{i}-r_{0,i}-\tau
^{2}\zeta _{i}\right) \right\Vert  \nonumber \\
&\leq & \left( \left\Vert r_{1,i}\right\Vert +\tau ^{2}\left\Vert
Lr_{2,i}\right\Vert +\left\Vert Lr_{3,i}\right\Vert +\left\Vert r_{4,i}-\tau
^{2}\zeta _{i}\right\Vert \right) \ .  \label{R2.1}
\end{eqnarray}

By formula (\ref{E1.35}), we get
\begin{eqnarray}
\left\Vert r_{1,i}\right\Vert &=&\left\Vert \left( S-L\right) \left[
u(t_{i})-u(t_{i-1})\right] \right\Vert  \nonumber \\
&=&\tau ^{2}\left\Vert \sum_{j=1}^{m}\eta _{j}^{-1}\left( I-S_{j}\right)
\left( \eta _{j}^{-1}A_{j}A^{-1}-I\right) AL\left[ u(t_{i})-u(t_{i-1})\right]
\right\Vert  \nonumber \\
&\leq &\tau ^{2}c_{4}\left\Vert A\left[ u(t_{i})-u(t_{i-1})\right]
\right\Vert ,  \label{zk2.16}
\end{eqnarray}

Obviously, for $Lr_{2,i}$ we have
\begin{equation}
\left\Vert Lr_{2,i}\right\Vert \leq \left\Vert A\left[ u(t_{i+1})-u(t_{i})%
\right] \right\Vert .  \label{zk2.17}
\end{equation}

From inequalities (\ref{zk2.16}) and (\ref{zk2.17}) , with account of
equation (\ref{2.1}), we get
\begin{eqnarray}
&&\left\Vert r_{1,i}\right\Vert +\tau ^{2}\left\Vert Lr_{2,i}\right\Vert
\nonumber \\
&\leq &\tau ^{2}c_{4}\left( \left\Vert u^{\prime \prime }\left( t_{i}\right)
-u^{\prime \prime }\left( t_{i-1}\right) \right\Vert +\left\Vert
f(t_{i})-f(t_{i-1})\right\Vert \right)  \nonumber \\
&&+\tau ^{2}\left( \left\Vert u^{\prime \prime }\left( t_{i+1}\right)
-u^{\prime \prime }\left( t_{i}\right) \right\Vert +\left\Vert
f(t_{i+1})-f(t_{i})\right\Vert \right) \ .  \label{zk2.19}
\end{eqnarray}

It is also obvious, that from (\ref{R4}) it follows
\begin{equation}
\left\Vert r_{4,i}-\tau ^{2}\zeta _{i}\right\Vert \leq \tau ^{2}\left\Vert
S_{1}\left[ M\left( u(t_{i})\right) -M\left( v_{i}\right) \right]
\right\Vert \leq a\tau ^{2}\left\Vert u(t_{i})-v_{i}\right\Vert .
\label{zk2.21}
\end{equation}

From the representation of $r_{3,i}$ it follows that
\begin{eqnarray}
\left\Vert Lr_{3,i}\right\Vert &\leq
&\int_{t_{i}}^{t_{i+1}}(t_{i+1}-t)\left\Vert u^{\prime \prime }\left(
t\right) -u^{\prime \prime }\left( t_{i}\right) \right\Vert dt  \nonumber \\
&&+\int_{t_{i-1}}^{t_{i}}(t-t_{i-1})\left\Vert u^{\prime \prime }\left(
t\right) -u^{\prime \prime }\left( t_{i}\right) \right\Vert dt  \nonumber \\
&\leq &\tau \int_{t_{i-1}}^{t_{i+1}}\left\Vert u^{\prime \prime }\left(
t\right) -u^{\prime \prime }\left( t_{i}\right) \right\Vert dt\ .
\label{zk2.23}
\end{eqnarray}

If we substitute inequalities (\ref{zk2.19}), (\ref{zk2.21}) and (\ref%
{zk2.23}) in (\ref{R2.1}), we get
\begin{eqnarray}
&&\left\Vert \left( R_{k-i}-R_{k-1-i}\right) \left( r_{i}-r_{0,i}-\tau
^{2}\zeta _{i}\right) \right\Vert  \nonumber \\
&\leq &\tau ^{2}c_{4}\left( \left\Vert u^{\prime \prime }\left(
t_{i}\right) -u^{\prime \prime }\left( t_{i-1}\right) \right\Vert
+\left\Vert f(t_{i})-f(t_{i-1})\right\Vert \right)  \nonumber \\
&&+\tau ^{2}\left( \left\Vert u^{\prime \prime }\left(
t_{i+1}\right) -u^{\prime \prime }\left( t_{i}\right) \right\Vert
+\left\Vert f(t_{i+1})-f(t_{i})\right\Vert \right)  \nonumber \\
&&+\tau \int_{t_{i-1}}^{t_{i+1}}\left\Vert u^{\prime \prime }\left(
t\right) -u^{\prime \prime }\left( t_{i}\right) \right\Vert dt+ac_{4}\tau
^{2}\left\Vert u(t_{i})-v_{i}\right\Vert \ .  \label{zk2.24}
\end{eqnarray}

According to estimate (\ref{B2}), we have
\begin{eqnarray}
&&\left\Vert R_{k-i}\left( r_{0,i}-r_{0,i-1}\right) \right\Vert  \nonumber \\
&=&\left\Vert R_{k-i}\left( I-B^{2}\right) ^{1/2}\left( I-B^{2}\right)
^{-1/2}\left( r_{0,i}-r_{0,i-1}\right) \right\Vert  \nonumber \\
&\leq &\left\Vert U_{k-i}\left( I-B^{2}\right) ^{1/2}\right\Vert \left\Vert
\left( I-S\right) ^{-1/2}\left( r_{0,i}-r_{0,i-1}\right) \right\Vert
\nonumber \\
&\leq &\left\Vert \left( I-S\right) ^{-1/2}\left( r_{0,i}-r_{0,i-1}\right)
\right\Vert .  \label{zk2.25}
\end{eqnarray}

Similarly we get
\begin{eqnarray}
&&\left\Vert R_{k-i}\left( \zeta _{i}-\zeta _{i-1}\right) \right\Vert
\nonumber \\
&\leq &\left\Vert \left( I-S\right) ^{-1/2}\left( \zeta _{i}-\zeta
_{i-1}\right) \right\Vert  \nonumber \\
&\leq &\left\Vert \left( I-S\right) ^{-1/2}(S_{1}-I)\left( \widetilde{f}%
(t_{i})-\widetilde{f}(t_{i-1})\right) \right\Vert  \nonumber \\
&&+\tau ^{2}\left\Vert \left( I-S\right) ^{-1/2}AL\left( \widetilde{f}%
(t_{i})-\widetilde{f}(t_{i-1})\right) \right\Vert .  \label{zk2.27}
\end{eqnarray}

If we substitute in (\ref{zk2.25}) representation of $r_{0,i}$ ($r_{0,i}=$ $%
\left( S-L\right) u(t_{i})$ ) and take into account estimate (\ref{E1.43}),
we get
\begin{equation}
\left\Vert R_{k-i}\left( r_{0,i}-r_{0,i-1}\right) \right\Vert \leq \tau
^{2}c_{1}\left\Vert A\left[ u(t_{i})-u(t_{i-1})\right] \right\Vert \ .
\label{zk2.28}
\end{equation}

Similarly we have
\begin{equation}
\left\Vert R_{k-1}r_{0,0}\right\Vert \leq \tau ^{2}c_{1}\left\Vert
Au(0)\right\Vert \ =\tau ^{2}c_{1}\left\Vert A\varphi _{0}\right\Vert \
,\quad \varphi _{0}\in D(A).  \label{zk2.281}
\end{equation}

From (\ref{zk2.28}), taking into account equation (\ref{2.1}), we have
\begin{eqnarray}
&&\left\Vert R_{k-i}\left( r_{0,i}-r_{0,i-1}\right) \right\Vert  \nonumber \\
&\leq &\tau ^{2}c_{1}\left( \left\Vert u^{\prime \prime }\left( t_{i}\right)
-u^{\prime \prime }\left( t_{i-1}\right) \right\Vert +\left\Vert
f(t_{i})-f(t_{i-1})\right\Vert \right) \ .  \label{zk2.29}
\end{eqnarray}

If in (\ref{zk2.27}) we take into account estimates (\ref{E1.59}), (\ref%
{E1.39}) and $\tau \left\Vert A^{1/2}L^{1/2}h\right\Vert \leq h$, we get
\begin{eqnarray}
&&\left\Vert R_{k-i}\left( \zeta _{i}-\zeta _{i-1}\right) \right\Vert
\nonumber \\
&\leq &\eta _{1}^{-1/2}\left\Vert \widetilde{f}(t_{i})-\widetilde{f}%
(t_{i-1})\right\Vert +\tau m\left\Vert A^{1/2}L\left( \widetilde{f}(t_{i})-%
\widetilde{f}(t_{i-1})\right) \right\Vert  \nonumber \\
&&+\tau c_{0}\left\Vert A^{1/2}L^{1/2}\left( \widetilde{f}(t_{i})-\widetilde{%
f}(t_{i-1})\right) \right\Vert \leq c_{6}\left\Vert \widetilde{f}(t_{i})-%
\widetilde{f}(t_{i-1})\right\Vert  \nonumber \\
&\leq &c_{6}\left( \left\Vert f(t_{i})-f(t_{i-1})\right\Vert +a\left\Vert
u(t_{i})-u(t_{i-1})\right\Vert \right) \ ,  \label{zk2.31}
\end{eqnarray}
where $c_{6}=\eta _{1}^{-1/2}+m+c_{0}$ .

Analogously to (\ref{zk2.31}), the inequality is true
\begin{equation}
\left\Vert R_{k-1}\zeta _{0}\right\Vert \leq c_{6}\left\Vert \widetilde{f}%
(t_{0})\right\Vert \leq c_{6}\left( \left\Vert f(0)\right\Vert +\left\Vert
M(\varphi _{0})\right\Vert \right) \ .  \label{zk2.32}
\end{equation}

\qquad Obviously, from inequalities (\ref{zk2.29}) and (\ref{zk2.31}) it
follows
\begin{eqnarray}
&&\left\Vert R_{k-i}\left[ \left( r_{0,i}-r_{0,i-1}\right) +\tau ^{2}\left(
\zeta _{i}-\zeta _{i-1}\right) \right] \right\Vert  \nonumber \\
&\leq &\tau ^{2}\left( c_{1}\left\Vert u^{\prime \prime }\left( t_{i}\right)
-u^{\prime \prime }\left( t_{i-1}\right) \right\Vert +ac_{6}\left\Vert
u(t_{i})-u(t_{i-1})\right\Vert \right)  \nonumber \\
&&+\tau ^{2}(c_{1}+c_{6})\left\Vert f(t_{i})-f(t_{i-1})\right\Vert \ .
\label{zk2.35}
\end{eqnarray}

If in (\ref{zk2.15}) we move to norms and take into account
inequalities (\ref{zk2.24}), (\ref{zk2.35}), (\ref{zk2.281}) and (\ref%
{zk2.32}), we get
\begin{eqnarray*}
\left\Vert \Phi _{k}\right\Vert &\leq &\tau ^{2}c_{4}\sum_{i=1}^{k}\left( \left\Vert u^{\prime \prime }\left(
t_{i}\right) -u^{\prime \prime }\left( t_{i-1}\right) \right\Vert
+\left\Vert f(t_{i})-f(t_{i-1})\right\Vert \right) \\
&&+\tau ^{2}\sum_{i=1}^{k}\left( \left\Vert u^{\prime \prime }\left(
t_{i+1}\right) -u^{\prime \prime }\left( t_{i}\right) \right\Vert
+\left\Vert f(t_{i+1})-f(t_{i})\right\Vert \right) \\
&&+\sum_{i=1}^{k}\left( \tau \int_{t_{i-1}}^{t_{i+1}}\left\Vert
u^{\prime \prime }\left( t\right) -u^{\prime \prime }\left( t_{i}\right)
\right\Vert dt+a\tau ^{2}\left\Vert z_{i}\right\Vert \right)
\end{eqnarray*}

\begin{eqnarray*}
&&+\tau ^{2}\sum_{i=1}^{k}\left( c_{1}\left\Vert u^{\prime \prime }\left(
t_{i}\right) -u^{\prime \prime }\left( t_{i-1}\right) \right\Vert
+ac_{6}\left\Vert u(t_{i})-u(t_{i-1})\right\Vert \right) \\
&&+\tau ^{2}(c_{1}+c_{6})\sum_{i=1}^{k}\left\Vert
f(t_{i})-f(t_{i-1})\right\Vert \\
&&+\tau ^{2}c_{1}\left\Vert A\varphi _{0}\right\Vert +\tau ^{2}c_{6}\left(
\left\Vert f(0)\right\Vert +\left\Vert M(\varphi _{0})\right\Vert \right) \ .
\end{eqnarray*}

Or the same
\begin{eqnarray}
\left\Vert \Phi _{k}\right\Vert &\leq &\tau ^{2}\sum_{i=1}^{k}\left[
c_{5}J_{i}\left( t_{i-1},u^{\prime \prime }\right) +J_{i}\left(
t_{i+1},u^{\prime \prime }\right) +ac_{6}J_{i}\left( t_{i-1},u\right) \right]
\nonumber \\
&&+\tau ^{2}\sum_{i=1}^{k}\left[ c_{7}J_{i}\left( t_{i-1},f\right) +J_{i}\left( t_{i-1},f\right) \right]  \nonumber \\
&&+\tau \sum_{i=1}^{k}\left( \int_{t_{i-1}}^{t_{i+1}}J_{i}\left(
t,u^{\prime \prime }\right) dt+\tau a\left\Vert z_{i}\right\Vert \right)
\nonumber \\
&&+\tau ^{2}\left[ c_{1}\left\Vert A\varphi _{0}\right\Vert +c_{6}\left(
\left\Vert f(0)\right\Vert +\left\Vert M(\varphi _{0})\right\Vert \right) %
\right] \ ,  \label{zk2.39}
\end{eqnarray}
where $c_{7}=c_{4}+c_{1}+c_{6}$ .

From (\ref{zk2.3}), taking into account (\ref{zk2.5}), (\ref{zk2.7})
and (\ref{zk2.39}), it follows (\ref{zkdel}).

\textbf{Result 2.} If functions $f(t)$ and $u^{\prime \prime }\left(
t\right) $ on the interval $[0,\ T]$ satisfies Holder condition with the
index $\lambda $ ($0<\lambda \leq 1$) , then

\begin{equation}
\left\Vert u^{\prime }\left( t_{k}\right) -\frac{v_{k+1}-v_{k}}{\tau }\right\Vert \leq c\tau ^{\lambda },\quad c=const>0.  \label{zk2.45}
\end{equation}

With account of (\ref{zkdel}), estimate (\ref{zk2.45}) follows from the following equality
\begin{eqnarray*}
&&u^{\prime }\left( t_{k}\right) -\frac{v_{k+1}-v_{k}}{\tau } \\
&=&\tau ^{-1}\int_{t_{k}}^{t_{k+1}}\left[ u^{\prime }\left(t_{k}\right) -u^{\prime }\left( t\right) \right] dt-\frac{z_{k+1}-z_{k}}{\tau }
\end{eqnarray*}

\section{Numerical Results}
The calculations were performed for the following problem.
\begin{eqnarray*}
&&\frac{\partial^{2}u}{\partial t^{2}} - \Delta u = sin(u) + f(x,y,t), \ \ (x,y,t) \in \Omega \times (0,T), \\
&&u(x,y,0)=\varphi_{0}(x,y), \ \ \ u'(x,y,0)=\varphi_{1}(x,y), \\
&&\Omega = (0,1) \times (0,1), \\
&&u_{\partial \Omega \times [0,T] } = 0.
\end{eqnarray*}

For this problem two different test case were calculated.

{\bf Test 1.}
\begin{eqnarray*}
&&u\left( x,y,t \right) = t^{7/2}\sin\left( 2{\pi}{x} \right)\cos\left( 2{\pi}{x} \right),\\
&&(x,y,t) \in [0,1]\times[0,1]\times[0,1],\\
&&\tau=h_{x}=h_{y}=0.05.
\end{eqnarray*}

\begin{figure}[H]
	\centering
	\includegraphics[width=0.9\textwidth]{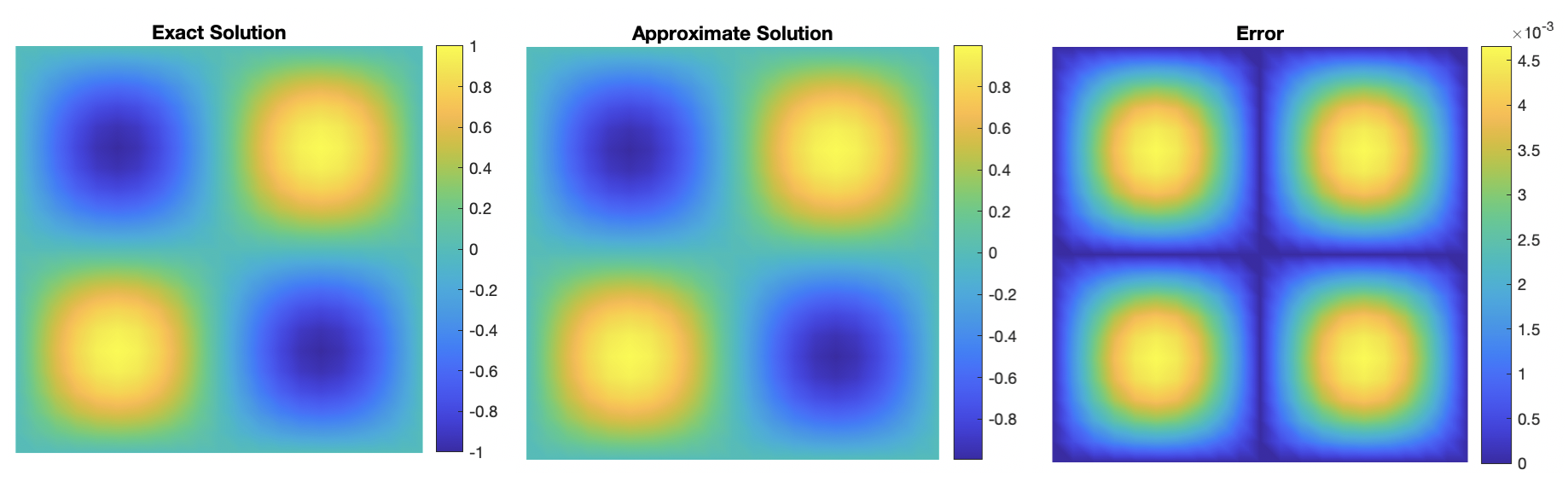}
	\caption{Exact and approximate solutions of Test 1.}
	\label{fig:figure1}
\end{figure}

Fig \ref{fig:figure1} shows exact solution, approximate solution and error. This test example is interesting as fourth derivative of the solution is discontinuous. As we see from Fig. \ref{fig:figure1} the numerical algorithm was able to resolve this problem with good accuracy. Maximum error equals to $0.0046$.

{\bf Test 2.}
\begin{eqnarray*}
&&u\left( x,y,t \right) = t^{7/2}\sin\left( 10{\pi}{x} \right)\cos\left( 10{\pi}{x} \right),\\
&&(x,y,t) \in [0,1]\times[0,1]\times[0,1],\\
&&\tau=h_{x}=h_{y}=0.01.
\end{eqnarray*}

\begin{figure}[H]
	\centering
	\includegraphics[width=0.9\textwidth]{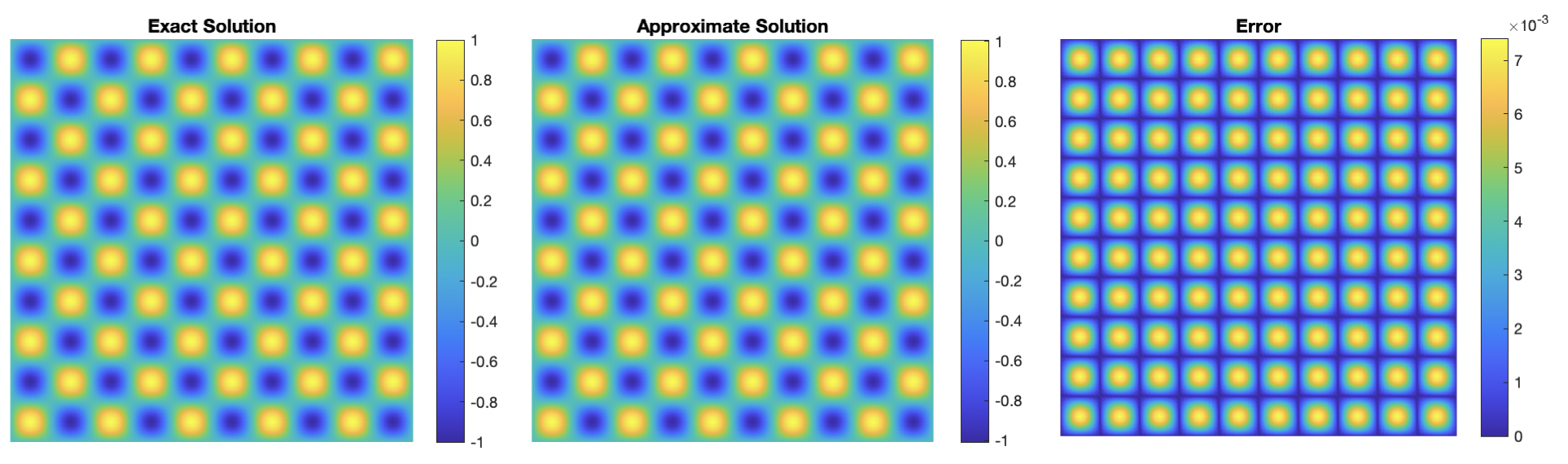}
	\caption{Exact and approximate solutions of Test 2.}
	\label{fig:figure2}
\end{figure}

Similarly, Fig \ref{fig:figure2} shows exact solution, approximate solution and error. This test example is interesting as fourth derivative of the solution is discontinuous and also has high oscillations, in the interval $\left( 0,1 \right)$ it changes the sign ten times. As we see from Fig. \ref{fig:figure2} the numerical algorithm was able to detect these oscillations and recovered original function with good accuracy. Maximum error for this case equals to $0.0074$.

\begin{figure}[H]
	\centering
	\includegraphics[width=0.9\linewidth]{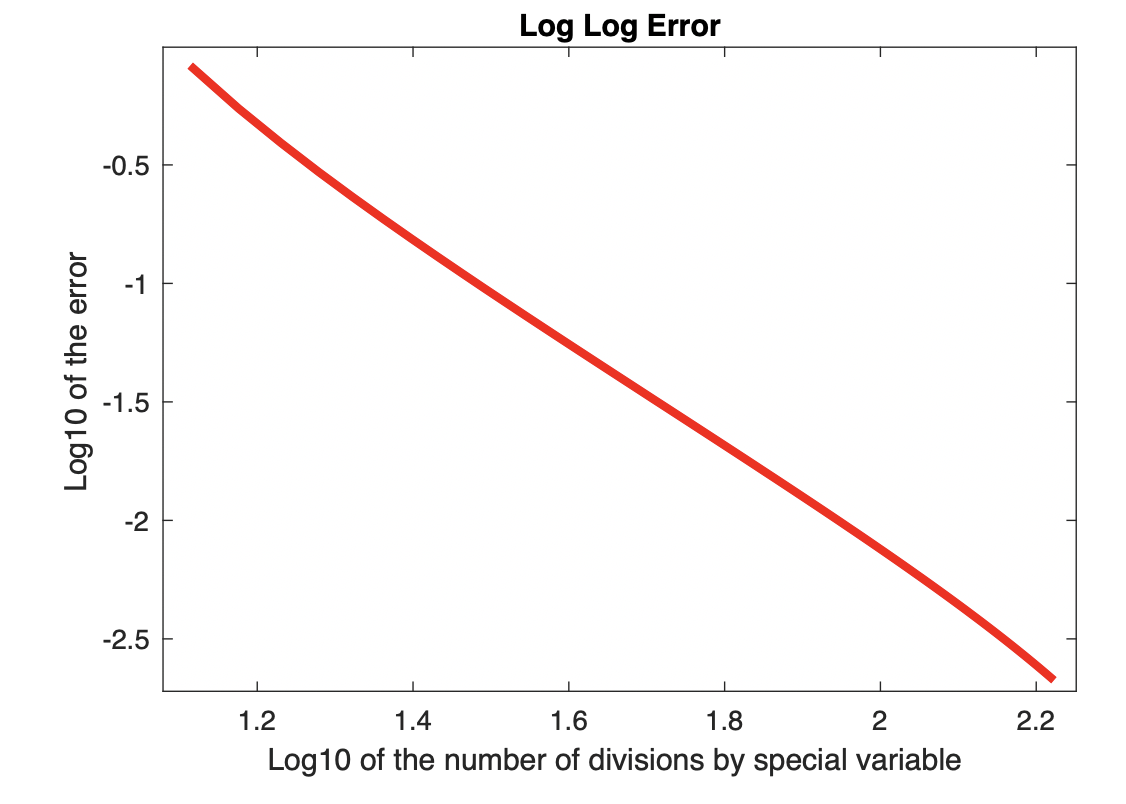}
	\caption{Dependence of logarithm of relative error on logarithm of number of divisions for Test 2.}
	\label{fig:figure3}
\end{figure}

On Fig. \ref{fig:figure3}, there is given a dependence of the logarithm of the relative error of the approximated solution on the logarithm of number of division by special variable. Aim of this figure is to find the convergence rate of the method by means of the
numerical experiment. If the method is second order of accuracy, then, the curve of the function (logarithm of the solution error) should approach to the line, the tangent of which equals two. On Fig. \ref{fig:figure3} it is clearly seen that, the curve approaches the line, the tangent of which equals to two, and this verifies the theoretical result proved in the article.

\def\printchapternonum{}
\bibliographystyle{plain}
\bibliography{refsource}

\section*{Authors' addresses:}

\begin{description}
	\item[{\Dikhaminjia}] \hfill \\ School of Business, Technology and Education, Ilia State University (ISU), Kakutsa Cholokashvili Ave 3/5, Tbilisi 0162, Georgia. \\ E-mail: \href{mailto:nana.dikhaminjia@iliauni.edu.ge}{\textbf{nana.dikhaminjia@iliauni.edu.ge}}
	\item[{\Rogava}] \hfill \\ Faculty of Exact and Natural Sciences, Ivane Javakhishvili Tbilisi State University (TSU), Ilia Vekua Institute of Applied Mathematics (VIAM), 2 University St., Tbilisi 0186, Georgia. \\ E-mail: \href{mailto:jemal.rogava@tsu.ge}{\textbf{jemal.rogava@tsu.ge}}
	\item[{\Tsiklauri}] \hfill \\ Missouri University of Science and Technology, Electromagnetic Compatibility Laboratory, 4000 Enterprise Drive, Rolla, MO 65409, USA. \\ E-mail: \href{mailto:tsiklaurim@mst.edu}{\textbf{tsiklaurim@mst.edu}}, \href{mailto:mtsiklauri@gmail.com}{\textbf{mtsiklauri@gmail.com}}
\end{description}

\end{document}